\documentclass{siamart1116}
\usepackage{algorithmic}
\usepackage{algorithm}

\usepackage{lipsum}
\usepackage{amsfonts}
\usepackage{graphicx}
\usepackage{epstopdf}
\usepackage{algorithmic}
\ifpdf
  \DeclareGraphicsExtensions{.eps,.pdf,.png,.jpg}
\else
  \DeclareGraphicsExtensions{.eps}
\fi

\numberwithin{theorem}{section}

\newcommand{\TheTitle}{Stabilization of BiCGSTAB by the generalized residual cutting method}
\newcommand{\TheAuthors}{T. Abe}

\headers{\TheTitle}{\TheAuthors}

\title{{\TheTitle}\thanks{Submitted to the editors DATE.}
}
\author{T. Abe\thanks{Kyoritsu Women's Junior College,
  Chiyoda-ku, Tokyo 101-8437, Japan
  }
    (\email{abe.toshihiko24@gmail.com})
}

\usepackage{amsopn}

\ifpdf
\hypersetup{
  pdftitle={\TheTitle},
  pdfauthor={\TheAuthors}
}
\fi

\newcommand{\eq}{\begin{equation}}
\newcommand{\qe}{\end{equation}}
\newcommand{\bgc}{\begin{center}}
\newcommand{\edc}{\end{center}}

\begin{document}

\maketitle

\begin{abstract}
The residual cutting (RC) method has been
proposed as an outer-inner loop iteration for efficiently solving
large and sparse linear systems of equations arising in solving numerically
problems of elliptic partial differential
equations. Then based on RC the generalized residual cutting (GRC)
method has been introduced, which can be applied to more general sparse linear systems problems.
In this paper, we
show that GRC can stabilize the BiCGSTAB, which is also an iterative algorithm
for solving large, sparse, and nonsymmetric linear systems,
and widely used in scientific computing and engineering simulations, due to its robustness.
BiCGSTAB converges faster and more smoothly than the original BiCG method,
by reducing irregular convergence behavior by stabilizing residuals.
However, it sometimes fails to converge due to stagnation or breakdown.
We attempt to enhance its robustness by further stabilizing it by GRC,
avoiding such failures.

\end{abstract}



\section{Introduction}
In \cite{grc2023}, we proposed a stable version of GRC with modified Gram-Schmidt algorithm.
In this report, based on GRC, we aim to construct a method which can avoid failures of
BiCGSTAB, such as stagnation or breakdown of the algorithm.
Thus, we put our emphasis on stability of calculation, rather than speeding up the convergence.

\subsection{Modified Gram-Schmidt based GRC}
Let us consider solving a linear system
\eq
{\bf H} u = b. 
\qe
The GRC (stabilized by MGS) is presented in Algorithm 1,
where $\psi^m$ is the solution of the inner loop iteration at the outer loop iteration $m$
and $\phi^m$
is a linear combination of vectors
$\psi^m$ and $\{ \phi^{m-k+1} \}, k=2, \cdots, j$.

\begin{algorithm}
  \caption{GRC (based on MGS)}
  \label{alg:mgsgrc}
  \begin{algorithmic}
    \STATE{1. Input: initial vectors ${u}^0$, ${r}^0 = {b} - {\bf H} {u}^0$}
    \STATE{2. For $m=0$ Until Convergence Do (outer loop)}
    \STATE{\hspace{10mm} 3. Compute approximate solution ${\psi}^m$ for ${\bf H} {\psi}^m = {r}^m$, \\
    \hspace{20mm} using ${\psi}^m = G({\phi}^{m-1},{\bf H},{r}^m)$ (inner loop)}
    \STATE{\hspace{10mm} 4. We apply Algorithm MGS to the set of vectors \\
    \hspace{20mm} $\{ {\bf H} \psi^m, \{ {\bf H} \phi^{m-k+1} \}, k=2, \cdots, j \}$ and compute \\
    \hspace{20mm} ${\bf H} \phi^m =$ LC$\{  {\bf H} \psi^m, \{ {\bf H} \phi^{m-k+1} \}, k=2, \cdots, j \}$ and \\
    \hspace{20mm} $\phi^m =$ LC$\{ \psi^m, \{ \phi^{m-k+1} \}, k=2, \cdots, j \}$
    }
    \STATE{\hspace{10mm} 5. $r^{m+1} = r^m - \alpha^m {\bf H} \phi^m$}
    where  $\alpha^m = (\psi^m, {\bf H} \phi^m) / ({\bf H} \phi^m, {\bf H} \phi^m)$
    \STATE{\hspace{10mm} 6. ${u}^{m+1} = {u}^{m} + \alpha^m \phi^m$}

    \STATE{End For}

    Output: $u^{m+1}$, $r^{m+1}$
  \end{algorithmic}
\end{algorithm}

\begin{algorithm}
  \caption{BiCGSTAB}
  \label{alg:bstab}
  \begin{algorithmic}
    \STATE{1. Input: initial vectors ${u}^0$, ${r}^0 = {b} - {\bf H} {u}^0$}
    \STATE{2. For $j=0$ Until $\frac{||r^{j+1}||}{||r^0||} < \theta$}
    \STATE{\hspace{10mm} 3. $\alpha_j = (r_j, {r_0}^{*}) / ({\bf H} p_j, {r_0}^*)$}
    \STATE{\hspace{10mm} 4. $s_j = r_j - \alpha_j{\bf H} p_j$}
    \STATE{\hspace{10mm} 5. $\omega_j = ({\bf H} s_j, s_j) / ({\bf H} s_j, {\bf H} s_j)$}
    \STATE{\hspace{10mm} 6. $x_{j+1} = x_j + \alpha_j p_j + \omega_j s_j$}
    \STATE{\hspace{10mm} 7. $r_{j+1} = s_j - \omega_j {\bf H} s_j$}
    \STATE{\hspace{10mm} 8. $\beta_j = \frac{(r_{j+1}, {r_0}^{*})}{(r_j, {r_0}^{*})}
                             \times \frac{\alpha_j}{\omega_j}$}
    \STATE{\hspace{10mm} 9. $p_{j+1} = r_{j+1} + \beta_j (p_j - \omega_j {\bf H} p_j)$}
    \STATE{End For}

    Output: $x^{j+1}$, $r^{j+1}$
  \end{algorithmic}
\end{algorithm}

\subsection{BiCGSTAB}
BiCGSTAB\cite{bicgstab} is derived by applying a stabilizing polynomial
to the CGS\cite{cgs}.
Its algorithm is described in Algorithm 2.

\section{GRC-BiCGSTAB}
Now we propose our method GRC-BiCGSTAB, which uses BiCGSTAB as
the inner loop in GRC.
Thus (some iterations of) Algorithm 2 is applied
to Algorithm 1 (line 3, inner loop) to compute
the approximate solution of
${\bf H} \psi^m = {r}^m$.

We set $\theta=0.5$ in in Algorithm 2 (line 2), which means
in the inner loop, BiCGSTAB iterates until the residual norm
decreases to the half of the initial one.
In the each step of the GRC outer loop, the solution is obtained as
the linear combination of the current BiCGSTAB solution output
and previous GRC solution vectors.
Thus the convergence of BiCGSTAB is further stabilized
by GRC (Please see Figure 1 for example.)



\section{Numerical experiments}
We validate the convergence of the proposed algorithm
and compare their performance with the original BiCGSTAB
and GRC,
by applying them to some important sparse linear systems
arising in engineering simulations problems.

\subsection{Experimental setup}
For numerical experiments, we used the BiCGSTAB function in the GNU Octave library.
We also implemented GRC-BiCGSTAB as an Octave function.
Thus, all calculations in our numerical tests use the same Octave routines such as matrix
times vector products, vector updates and dotproduct operations.
The stopping
criterion for convergence was that the Euclidean norm of the residual error was less than
$10^{-12}$ for both of BiCGSTAB and GRC.

We chose $j=5$ in Algorithm 1 for the GRC-BiCGSTAB, which was found out to be large enough
by experiments in \cite{grc2023}.
As already mentioned in the Section 2, we set $\theta=0.5$ for the BiCGSTAB
as the inner loop function.
Also, we set $\theta=10^{-12}$ for the original BiCGSSTAB, to meet with the stopping criterion.

\subsection{Experimental results}
The following test problems were used in the experiments.

\subsubsection{Problem 1}
A coefficient matrix generated by discretizing the partial differential
equation \cite{n11}
\[
u_{xx} + u_{yy} + u_{zz} + 1000u_x = F
\]
where $F(x,y,z) = \exp(x y z)  \sin(\pi x)  \sin(\pi y)  sin(\pi z)$.
The matrix of this problem is nonsymmetric and non positive definite of dimension
$n = 125$.

The result fo Problem 1 is shown in Figure 1.
The horizontal axis shows the number of iterations in BiCGSTAB.
For GRC-BiCGSTAB, it also shows the accumulated number of iterations
of the BiCGSTAB as the inner loop function, so that we can directly
compare the number of iterations in the same axis.
We can see the the GRC effectively stabilizes the fluctuation
of each inner loop calculation of BiCGSTAB.
In the each step of the GRC outer loop, the solution is obtained as
the linear combination of the current BiCGSTAB solution output
and previous GRC solution vectors.
Thus the convergence of BiCGSTAB is further stabilized
by GRC.

\begin{figure}
\bgc
\caption{Residual norms of BiCGSTAB and GRC-BiCGSTAB.
'BiCGSTAB' (blue) shows the original BiCGSTAB.
'GRC-BiCGSTAB' (green)  shows the residual norm stabilized by GRC.
'GRC-BiCGSTAB (inner loop)' (purple) shows the residual norm of BiCGSTAB
of the inner loop in GRC, where each residual at the beginning of the inner loop
is the one from the GRC outer loop.
}
\edc
\includegraphics[scale=0.58]{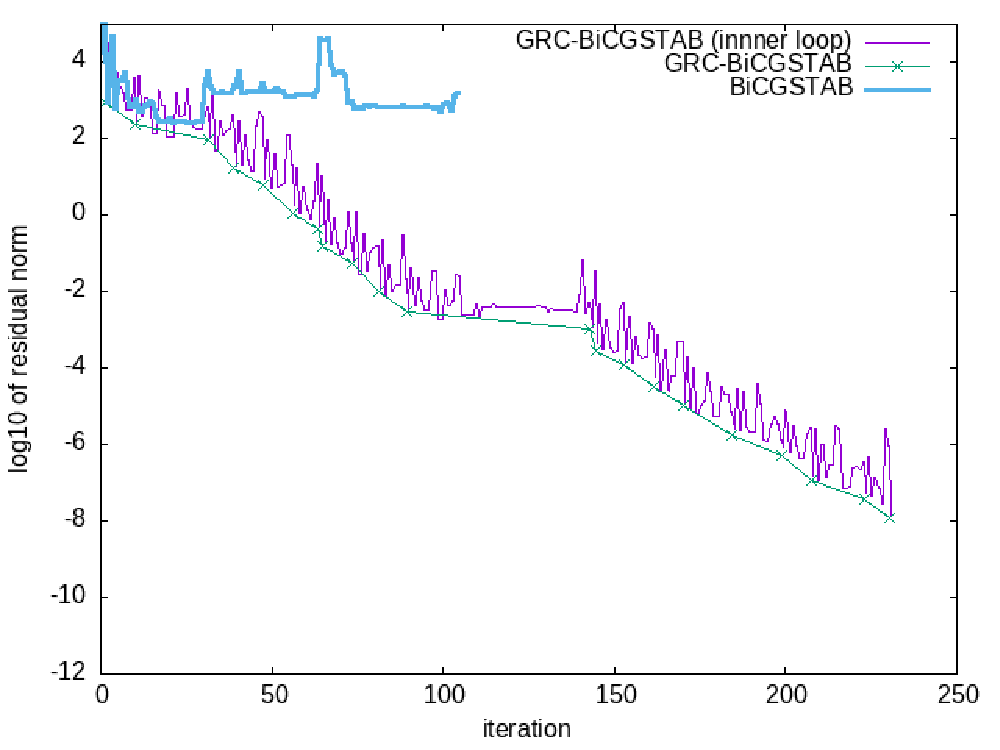}

Result with Problem 1. BiCGSTAB (blue) stopped due to breakdown at iteration 104,
while GRC-BiCGSTAB converged.

\end{figure}

\subsubsection{Other test problems}
Other problems are chosen from the widely used
database of sparse matrix collection \cite{n13}.
For these problems, the right hand side vector for the linear system (in eqn (1)) is selected so that the
solution vector will be the vector with all the elements
being unity (i.e. $u = [1, \cdots, 1]^T$) i.e. b is computed the product of the matrix and the solution
vector, whose elements are all unity.

\begin{figure}
\bgc\caption{Result for other test problems}\edc
\includegraphics[scale=0.5]{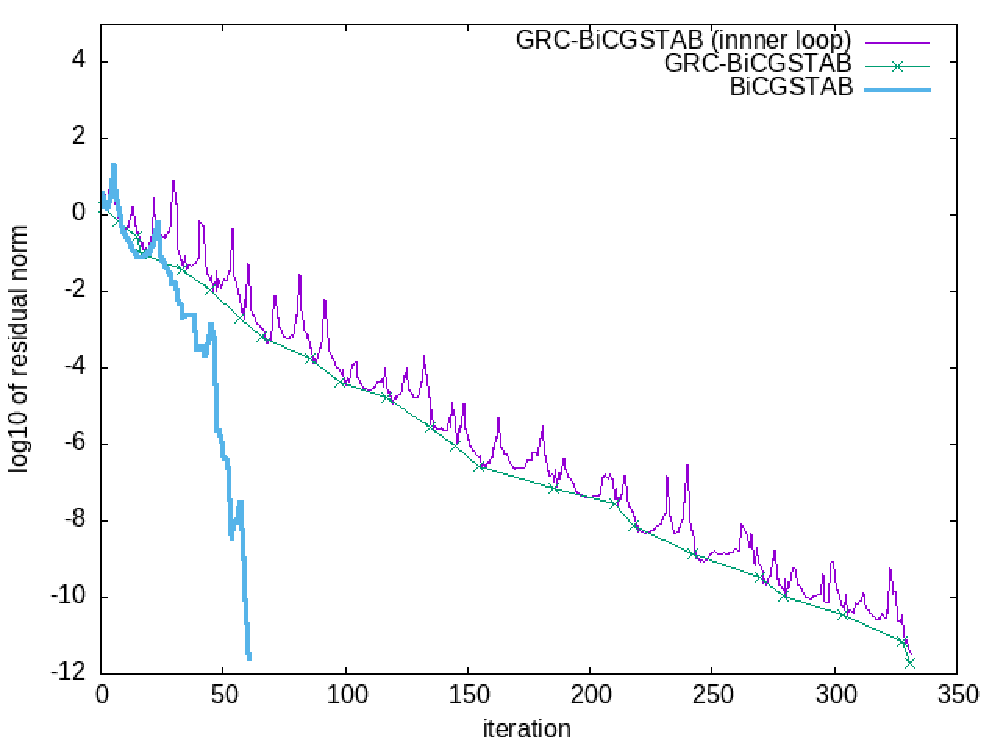}

(b) Result with bfwa62. In this case, BiCGSTAB converges significantly faster.

\includegraphics[scale=0.5]{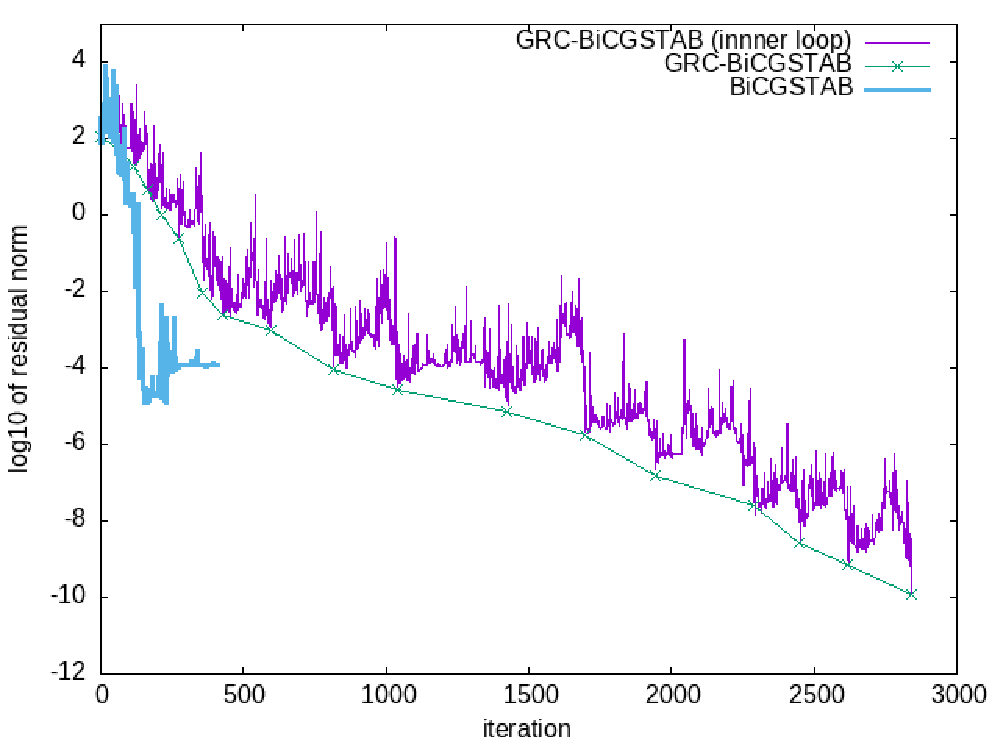}

(c) Result with olm100. BiCGSTAB stopped due to breakdown at the iteration 410.

\includegraphics[scale=0.5]{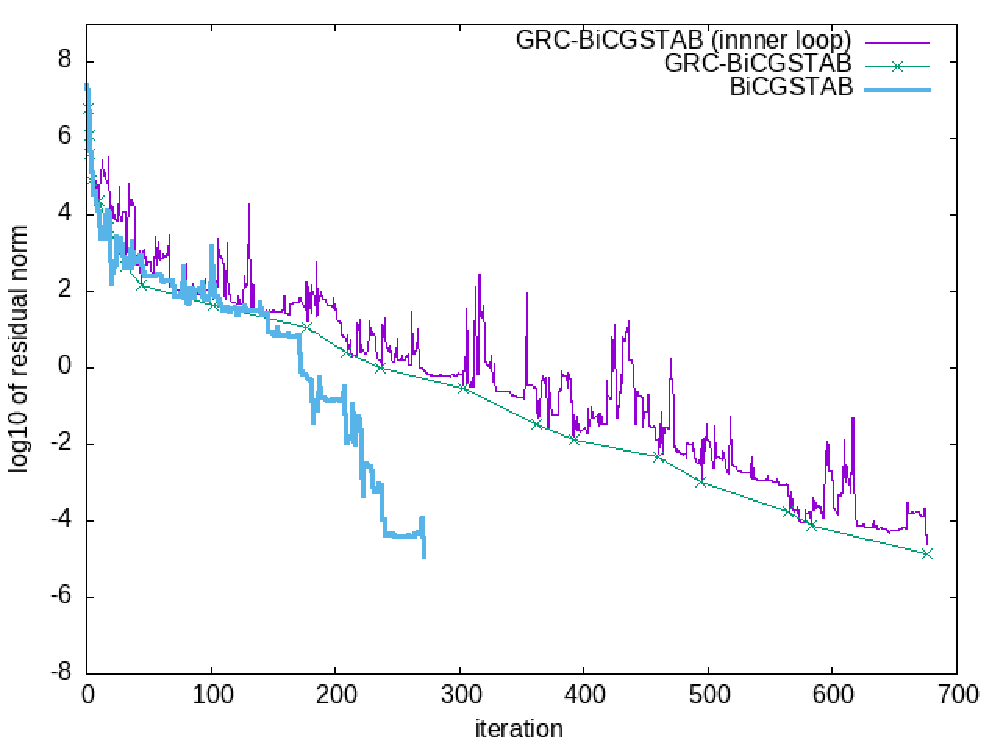}

(d) Result with pores\_1. BiCGSTAB converges significantly faster.

\end{figure}

\begin{figure}
\bgc\caption{Result for other test problems}\edc
\includegraphics[scale=0.5]{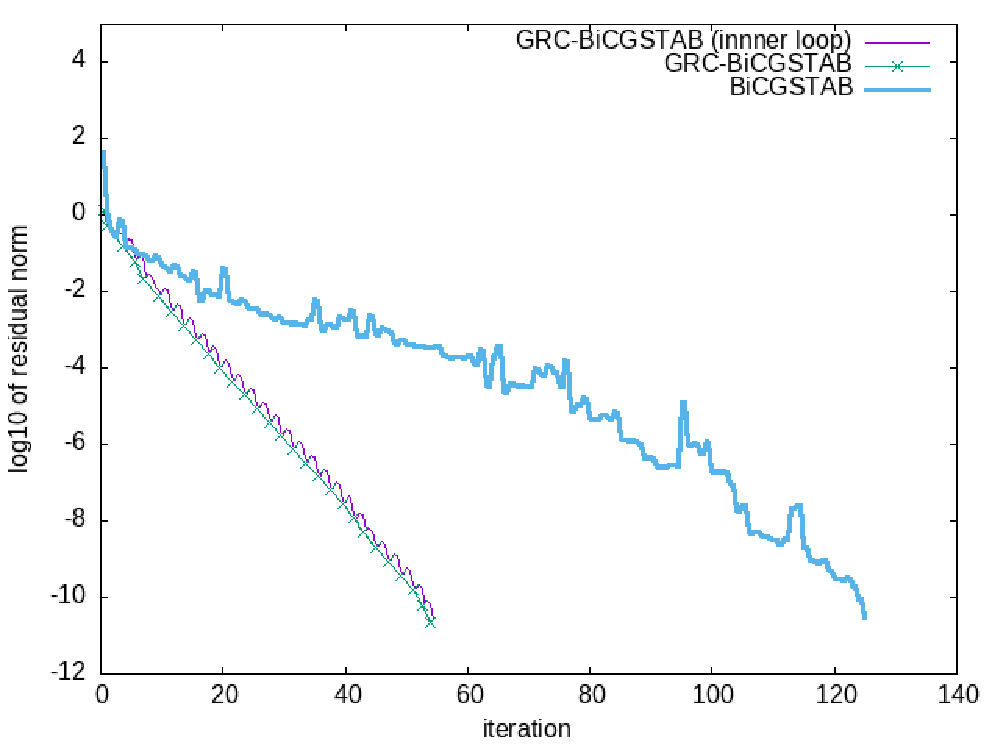}

(e) Result with toeplitz100\_1.4.

\includegraphics[scale=0.5]{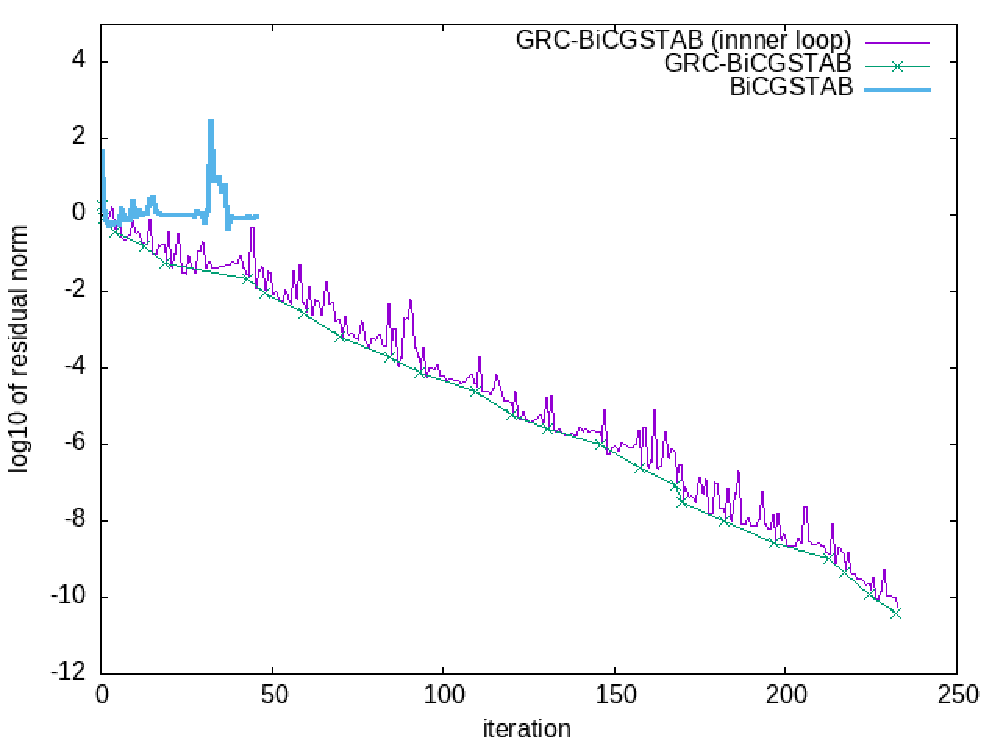}

(f) Result with toeplitz100\_2.0. BiCGSTAB stopped due to breakdown at the iteration 45.

\includegraphics[scale=0.5]{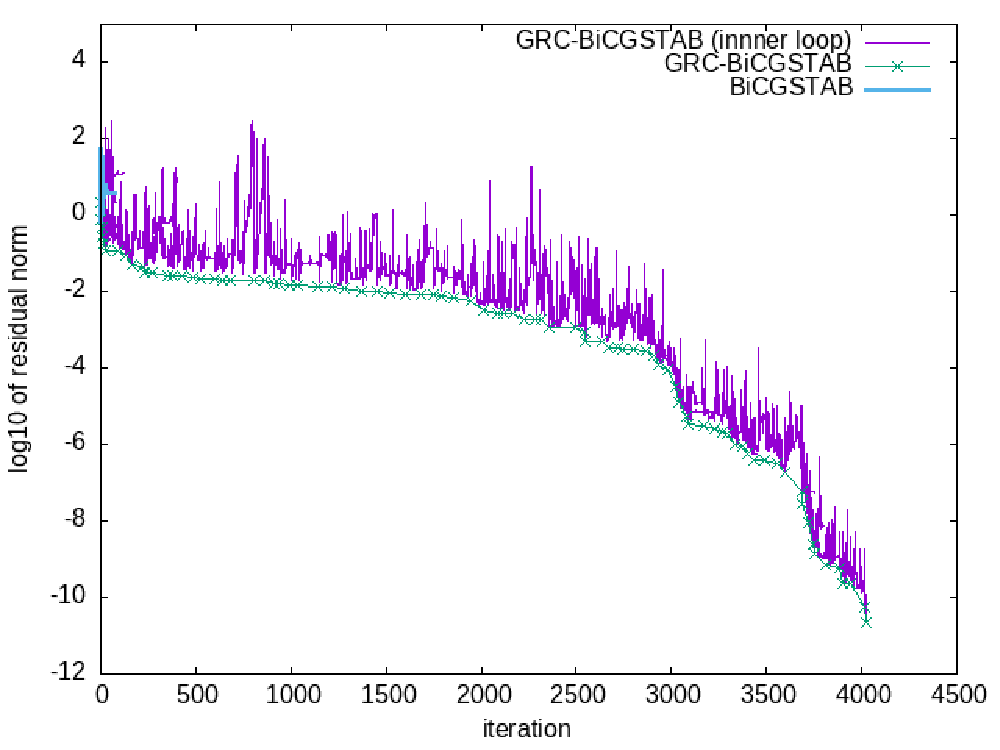}

(g) Result with toeplitz100\_2.3. BiCGSTAB stopped due to breakdown at the iteration 72.

\end{figure}

\section{Discussion and conclusion}
As shown in the experiments, BiCGSTAB converges faster
for the problems for which it converges,
while GRC-BiCGSTAB is more robust, so that it can converge
for the problems for which BiCGSTAB fails to converge.
As a whole, we consider GRC-BiCGSTAB is considered to be more suitable
when robust solution is needed.


\newpage

\begin{thebibliography}{9}

\bibitem{bicgstab}{\sc Van der Vorst, H. A.},
{Bi-CGSTAB: A Fast and Smoothly Converging Variant of Bi-CG for the Solution of Nonsymmetric Linear Systems}.
SIAM J. Sci. Stat. Comput. 13 (2) 1992: 631–644. 

\bibitem{cgs}{\sc Peter Sonneveld},
{CGS, A Fast Lanczos-Type Solver for Nonsymmetric Linear systems}.
SIAM Journal on Scientific and Statistical Computing. 10 (1) 1989: 36–52.

\bibitem{grc2023}{\sc Abe T, Chronopoulos AT},
{The generalized residual cutting method and its convergence characteristics}.
Numerical linear algebra with applications 2023; 30 (6): e2517.

\bibitem{n7}{\sc Saad Y},
{ Iterative methods for sparse linear systems}.
Society for Industrial and Applied Mathematics; 2003.

 \bibitem{n1} {\sc Tamura A, Kikuchi K, Takahashi},
{ Residual cutting method for elliptic boundary value problems: application to Poisson's equation}. J. Comp. Phys. 1997; 137: pp. 247-264.

\bibitem{n2} {\sc Abe T, Sekine Y, Sato F},
{ Solving a coupled perturbed equation by the residual cutting method}.
J. Chem. Phys. 2013; 557: 176-181.

\bibitem{n9} {\sc Abe T, Sekine Y, Kikuchi K},
{ Generalization of the residual cutting method based on the Krylov subspace}.
AIP Conf. Proc. 2016; 1738.
\bibitem{n13} {\sc Davis TA, Hu Y},
{ The University of Florida sparse matrix collection}.
ACM Transactions on Mathematical Software (TOMS) 2011; 38(1): 1-25.

\bibitem{n11} {\sc Sonneveld P, Gijzen MB},
{ IDR(s): A Family of simple and fast algorithms for solving large nonsymmetric systems of linear equations}.
SIAM Journal of Scientific Computing 2008; 31(2): 1035-1062.
\end{thebibliography}
\end{document}